\documentclass{amsart}
\usepackage{amsmath,amssymb,amscd}

\newcommand{\NI}{{\noindent}}
\newcommand{\QED}{\hfill$\Box$\medskip}
\newcommand{\ep}{\epsilon}

\newtheorem{theorem}{Theorem}[section]

\newtheorem{thm}[theorem]{Theorem}
\newtheorem{definition}[theorem]{Definition}

\newtheorem{rmk}[theorem]{Remark}
\newtheorem{lemma}[theorem]{Lemma}
\newtheorem{question}[theorem]{Question}
\newtheorem{guess}[theorem]{Conjecture}

\begin{document}

\title[The bounded isometry conjecture]{The bounded isometry conjecture for the Kodaira-Thurston
manifold and 4-Torus}
\author{Zhigang Han}
\address{Department of Mathematics and Statistics, University of Massachusetts,
Amherst, MA 01003-9305, USA} \email{han@math.umass.edu}
\date{May 02, 2007}

\begin{abstract}
The purpose of this note is to study the bounded isometry conjecture
proposed by Lalonde and Polterovich \cite{LP}. In particular, we
show that the conjecture holds for the Kodaira-Thurston manifold
with the standard symplectic form and for the 4-torus with all
linear symplectic forms.
\end{abstract}

\maketitle

\section{Introduction and main results}\label{introduction}
\smallskip

Let $(M, \omega)$ be a closed symplectic manifold. There is a
natural bi-invariant norm, called the Hofer norm $\rho$, defined on
the Hamiltonian diffeomorphism group ${\rm Ham}(M,\omega)$. That is,
$\rho(f)$ is the Hofer distance between the identity map $id$ and
$f$ for all $f \in {\rm Ham}(M,\omega)$, see Section \ref{hofernorm}
for details. Lalonde and Polterovich \cite{LP} have studied the full
symplectomorphism group ${\rm Symp}(M,\omega)$ within the framework
of Hofer's geometry. We first recall the notion of bounded and
unbounded symplectomorphisms. Namely, for each $\phi \in {\rm
Symp}(M,\omega)$, define
$$r(\phi):={\rm sup}\, \{\,\rho([\phi,f]) \mid f \in {\rm Ham}(M,\omega)\},$$
where $[\phi,f]:=\phi f \phi^{-1} f^{-1}$ is the commutator of
$\phi$ and $f$.

\smallskip

\begin{definition}\label{bounded}
An element $\phi \in {\rm Symp}(M,\omega)$ is bounded if
$r(\phi)<\infty$, and is unbounded if $r(\phi)= \infty$.
\end{definition}

Denote by  $BI_0(M,\omega)$ the set of all bounded elements in the
identity component ${\rm Symp}_0(M,\omega)$ of ${\rm
Symp}(M,\omega)$. Since $\rho$ is bi-invariant, it follows from the
inequality $\rho([\phi,f]) \leqslant 2 \rho(\phi)$ that ${\rm
Ham}(M,\omega)$ is a subgroup of $BI_0(M,\omega)$. The converse is
the following conjecture in \cite{LP}.

\smallskip

\begin{guess}[Bounded isometry conjecture]\label{bic} For all symplectic manifolds $(M,\omega)$,
$BI_0(M,\omega)={\rm Ham}(M,\omega)$.
\end{guess}

This conjecture was proved in \cite{LP} for closed surfaces with
area form and for arbitrary products of closed surfaces of genus
greater than 0 with product symplectic form; Lalonde and Pestieau
\cite{LPe} confirmed it for product symplectic manifolds $M = N
\times W$ with $N$ being any product of closed surfaces and $W$
being any closed symplectic manifold of first real Betti number
equal to zero.  In this note, we give a positive answer to this
conjecture for the Kodaira-Thurston manifold with the standard
symplectic form and for the 4-torus with all linear symplectic
forms.

\smallskip

\begin{thm}\label{bicforkt} The bounded isometry conjecture holds
for the Kodaira-Thurston manifold $M$ with the standard symplectic
form $\omega$.
\end{thm}

\begin{thm}\label{bicfortorus} The bounded isometry conjecture holds for the $4$-torus $(\mathbb
T^4, \omega)$ with any linear symplectic form $\omega:=\sum_{i <
j}\, a_{ij}\, dx_i \wedge dx_j$.
\end{thm}

\NI \textbf{Organization of the paper:} We begin with some
preparations in Section \ref{fluxsubgroup}-\ref{admissiblelift}.
Then  we prove Theorem \ref{bicforkt} in Section \ref{proofbicforkt}
and Theorem \ref{bicfortorus} in Section \ref{proofbicfortorus}. We
study the same conjecture for the Kodaira-Thurston manifold with all
linear symplectic forms in Section \ref{sec8}. While we are unable
to prove the conjecture in this case, some partial results are
provided and the difficulties are discussed.

\bigskip

\NI \textbf{Acknowledgements:} This work is part of the author's
Ph.D. thesis, being carried out under the supervision of Dusa McDuff
at Stony Brook University. The author wishes to thank her for her
great guidance and continual support. He also thanks Xiaojun Chen,
Basak Gurel, Leonid Polterovich, Guanyu Shi, Yujen Shu, Dennis
Sullivan and Aleksey Zinger for useful comments and discussions.

\bigskip

%%%%%%%%%%%%%%%%%%%%%%%%%%%%%%%%%%%%%%%%%%%%%%%%%%%%%%%%%%%%%%%%%%%%%%%
%%%%%%%%%%%%%%%%%%%%%%%%%%%%%%%%%%%%%%%%%%%%%%%%%%%%%%%%%%%%%%%%%%

\section{The flux subgroup}\label{fluxsubgroup}

\smallskip

The flux homomorphism is best defined on the universal cover
$\widetilde{\rm Symp}_0(M,\omega)$ of ${\rm Symp}_0(M,\omega)$,
$${\rm flux} : \widetilde{\rm Symp}_0(M,\omega) \to H^1(M, \mathbb R).$$
Let $\{\phi_t\} \in \widetilde{\rm Symp}_0(M,\omega)$, i.e. $\phi_t$
is a smooth isotopy in ${\rm Symp}_0(M,\omega)$. There exists a
unique family of vector fields $X_t$ which generates the flow
$\phi_t$, i.e.
$$\frac{d}{dt} \phi_t=X_t \circ \phi_t.$$
Define
$${\rm flux}(\{\phi_t\}):=\int_0^1 \iota(X_t)\omega \,dt.$$
In particular, if $\{\phi_t\}$ is the flow of the time-independent
symplecitc vector field $X$ on the time interval $0 \leqslant t
\leqslant 1$, then
\begin{equation} \label{eq:1}
\begin{split}
{\rm flux}(\{\phi_t\})=\iota(X)\omega.
 \end{split}
 \end{equation}
This fact will often be used in later calculations.

The flux subgroup $\Gamma:=\Gamma_\omega$ is the image
$${\rm flux}(\pi_1({\rm Symp}_0(M,\omega)) \subset H^1(M, \mathbb R)$$ of the
fundamental group of ${\rm Symp}_0(M,\omega)$ under the flux
homomorphism. Thus there is an induced map from ${\rm
Symp}_0(M,\omega)$, still denoted by ${\rm flux}$,
$${\rm flux} : {\rm Symp}_0(M,\omega) \to H^1(M, \mathbb
R)/\Gamma.$$ It is well known that this map is surjective, and its
kernel is equal to ${\rm Ham}(M,\omega)$. In other words, we have
the following exact sequence of groups
$$\begin{CD} 0 @>>> {\rm Ham}(M,\omega) @>>> {\rm Symp}_0(M,\omega)
@>{\rm flux}>>H^1(M, \mathbb R)/\Gamma @>>> 0. \end{CD}$$ We refer
to \cite{MS1} Chapter 10 for more details.

\smallskip

Since whether or not the flux is equal to 0 distinguishes a
Hamiltonian diffeomorphism from a nonHamiltonian symplectomorphism,
one main step in our applications is to understand the flux subgroup
$\Gamma$.

For this, we denote as in \cite{Ke} by $C(M)$ the space of
continuous maps from $M$ to $M$ with the compact open topology.
Given $p \in M$, we define the evaluation map $ev_c : C(M) \to M$ by
$ev_c(f)=f(p)$. Denote by $ev_s$ the restriction of $ev_c$ to ${\rm
Symp}_0(M,\omega)$. We will use the same notation for the induced
maps on the fundamental groups. By $\widetilde{ev}_s$ we denote the
homomorphism from $\pi_1({\rm Symp}_0(M,\omega))$ to $H_1(M, \mathbb
Z)$, which is the composition of $ev_s$ with the Hurewitz map from
$\pi_1(M)$ to $H_1(M,\mathbb Z)$.

The following commutative diagram due to Lalonde, McDuff and
Polterovich \cite{LMP2} plays a crucial role in the calculation of
the flux subgroup $\Gamma$.

\begin{lemma} [LMP] \label{commute}
Let $(M, \omega)$ be a closed symplectic manifold of dimension $2n$.
Then the following diagram commutes.
$$\begin{CD}
  \pi_1({\rm Symp}_0(M,\omega)) @>\widetilde{ev}_s>>H_1(M, \mathbb Z)    @>{\rm PD}>>H^{2n-1}(M, \mathbb Z) \\
  @VidVV             &&                      @VV\centerdot(n-1)!{\rm vol}(M)V\\
  \pi_1({\rm Symp}_0(M,\omega)) @>{\rm flux}>>H^1(M, \mathbb R) @>\wedge[\omega]^{n-1}>>  H^{2n-1}(M, \mathbb R). \end{CD}$$
\end{lemma}\QED

\bigskip

%%%%%%%%%%%%%%%%%%%%%%%%%%%%%%%%%%%%%%%%%%%%%%%%%%%%%%%%%%%%
%%%%%%%%%%%%%%%%%%%%%%%%%%%%%%%%%%%%%%%%%%%%%%%%%%%%%%%%%%%%%%%
\section{The Kodaira-Thurston manifold}\label{ktmanifold}

\smallskip

Let $G$ be the group $(\mathbb{Z}^4, \cdot)$ where $$(m_1, n_1, k_1,
\ell_1) \cdot (m_2, n_2, k_2, \ell_2)=(m_1+m_2, n_1+n_2, k_1+k_2+m_1
\ell_2, \ell_1+\ell_2).$$ $G$ acts on $\mathbb R^4$ via
$$G \to \mbox{Diff}(\mathbb R^4) : (m, n, k, \ell) \mapsto \rho_{m n k
\ell}$$ where
$$\rho_{m n k \ell}(s, t, x, y)=(s+m, t+n, x+k+m y, y+\ell).$$

Note that $\rho_{m n k \ell}$ preserves the symplectic form
$\omega=ds \wedge dt+dx\wedge dy$ on $\mathbb R^4$. Hence the
quotient $(M:=\mathbb R^4/G, \omega)$ is a closed symplectic
manifold, known as the Kodaira-Thurston manifold, see \cite{Th}. It
was the first known example of a closed symplectic manifold which
admits no k\"{a}hler structure, since its first betti number
$b_1=3$, see \cite{MS1} Example 3.8.

The manifold $M=\mathbb R^4/G$ can also be described as a torus
bundle over a torus, that is $M=\mathbb R^2 \times_{\mathbb{Z}^2}
\mathbb T^2$. Here $\mathbb{Z}^2$ acts on $\mathbb R^2$ in the usual
way, and it acts on $\mathbb T^2$ via
$$(m, n) \to A_{m n} : \left(
\begin{matrix}x\\y \end{matrix} \right)
 \mapsto \left( \begin{matrix}1&m\\0&1
\end{matrix} \right) \left( \begin{matrix}x\\y
\end{matrix} \right).$$ Therefore $M=\mathbb R \times S^1 \times
\mathbb T^2/\sim$, where
$$(s, t, x, y) \sim (s+1, t, x+y, y).$$

\smallskip

Our first task is to understand the flux subgroup $\Gamma$ of the
Kodaira-Thurston manifold described above. In particular, we have

\begin{thm}\label{flux} The flux subgroup $\Gamma \subset H^1(M,\mathbb R)$ of
the Kodaira-Thurston manifold with the standard symplectic form
$\omega=ds \wedge dt+dx\wedge dy$ has rank $2$ over $\mathbb{Z}$.
Namely, $\Gamma=\mathbb{Z}\langle ds, dy\rangle $.
\end{thm}

To prove Theorem \ref{flux}, we need the following result on the
cohomology groups of the Kodaira-Thurston manifold.

\begin{lemma}\label{cohomology}
The cohomology groups of the Kodaira-Thurston manifold $M$ described
above are as follows: $H^1(M, \mathbb R)$ is of rank $3$, generated
by $ds, \, dt$ and $dy$; $H^2(M, \mathbb R)$ is of rank $4$,
generated by $\gamma \wedge ds,\, \gamma \wedge dy,\, ds \wedge dt$
and $dy \wedge dt$; and $H^3(M, \mathbb R)$ is of rank $3$,
generated by $\gamma \wedge dy \wedge dt,\, \gamma \wedge dy \wedge
ds$ and $\gamma \wedge ds \wedge dt$, where $\gamma=dx-s dy$.
\end{lemma}

\begin{proof}
This follows from an easy calculation.
\end{proof}

\smallskip

\NI {\textbf{Proof of Theorem \ref{flux}.}} We use the commutative
diagram in Lemma \ref{commute}. For manifolds of dimension 4, the
diagram reads as
$$\begin{CD}
  \pi_1({\rm Symp}_0(M,\omega)) @>\widetilde{ev}_s>>H_1(M, \mathbb Z)    @>{\rm PD}>>H^3(M, \mathbb Z) \\
  @VidVV             &&                      @VV\centerdot {\rm vol}(M)V\\
  \pi_1({\rm Symp}_0(M,\omega)) @>{\rm flux}>>H^1(M, \mathbb R) @>\wedge[\omega]>> H^3(M, \mathbb R).
\end{CD}$$

\smallskip

Denote by $C_0(M)$ the identity component of $C(M)$. It was proved
in Gottlieb \cite{Go} (Theorem III.2) that for all aspherical
manifolds $M$,
$$ev_c : \pi_1(C_0(M))\cong Z(\pi_1(M))$$ is a group isomorphism,
where $Z(\pi_1(M))$ stands for the center of $\pi_1(M)$. For the
Kodaira-Thurston manifold $M=\mathbb R^4/G$, we have $\pi_1(M)=G$.
It is easy to check that
$Z(\pi_1(M))=\mathbb{Z}\langle\frac{\partial}{\partial t},
\frac{\partial}{\partial x}\rangle$, and the commutator group $[\,
\pi_1(M), \pi_1(M)]=\mathbb{Z}\langle\frac{\partial}{\partial
x}\rangle $, see Example 3.8 in \cite{MS1}. Thus the image of
$\widetilde{ev}_s$ in $H_1(M, \mathbb Z)$ is contained in
$$Z(\pi_1(M))/ [\, \pi_1(M),
\pi_1(M)]=\mathbb{Z}\langle\frac{\partial}{\partial t},
\frac{\partial}{\partial x}\rangle /
\mathbb{Z}\langle\frac{\partial}{\partial x}\rangle
=\mathbb{Z}\langle\frac{\partial}{\partial t}\rangle.$$

Note that $PD(\frac{\partial}{\partial t})=-dx \wedge dy \wedge
ds=-\gamma \wedge dy \wedge ds$, where $\gamma=dx-sdy$. Now look at
the map $\wedge \omega : H^1(M, \mathbb R) \to H^3(M, \mathbb R)$,
$$ds \mapsto ds \wedge \omega=\gamma \wedge dy \wedge ds \ne 0,$$
$$dt \mapsto dt \wedge \omega=\gamma \wedge dy \wedge dt \ne 0,$$
$$dy \mapsto dy \wedge \omega=dy \wedge ds \wedge dt=0.$$

\NI Here we have used the fact that the 3-form $dy \wedge ds \wedge
dt = d (\gamma \wedge dt)$ is exact, so it vanishes on the
cohomology level. Since ${\rm vol}(M)=1$, we conclude from the above
commutative diagram that the flux subgroup $\Gamma \subset H^1(M,
\mathbb R)$ is contained in $\mathbb{Z}\langle ds, dy\rangle $. An
explicit construction shows that $\Gamma$ is actually equal to
$\mathbb{Z}\langle ds, dy\rangle $. Namely, we take two elements
$\{\phi_\theta\}$ and $\{\psi_\theta\}$ in $\pi_1({\rm
Symp}_0(M,\omega))$ such that
$$\phi_\theta(s, t, x, y)=(s, t-\theta, x, y), 0 \leqslant\theta \leqslant1, $$
$$\psi_\theta(s, t, x, y)=(s, t, x+\theta, y), 0 \leqslant\theta \leqslant1.$$
Using (\ref{eq:1}) in Section \ref{fluxsubgroup}, one can show that
${\rm flux}(\{\phi_\theta \})=ds$ and ${\rm flux}(\{\psi_\theta
\})=dy$. This completes the proof of Theorem \ref{flux}.\QED

\bigskip

%%%%%%%%%%%%%%%%%%%%%%%%%%%%%%%%%%%%%%%%%%%%%%%%%%%%%%%%%%%%%%%
%%%%%%%%%%%%%%%%%%%%%%%%%%%%%%%%%%%%%%%%%%%%%%%%%%%%%%%%%%%%%%%%
\section{The Hofer norm}\label{hofernorm}

\smallskip

Let $(M,\omega)$ be a closed symplectic manifold of dimension $2n$.
Denote by $\mathcal{A}$ the space of all normalized smooth functions
on $M$ with respect to the volume form $\omega^n$, i.e.
$$\mathcal{A} := \{ F \in C^\infty (M) \mid \int_M F \,\omega^n
=0\}.$$ It is well known that $\mathcal{A}$ can be identified with
the space of Hamiltonian vector fields, which is the Lie
algebra\footnote{As a vector space, the Lie algebra is by definition
the tangent space to the Lie group at the identity. The tangent
spaces to the Lie group at other points are identified with the Lie
algebra with the help of right shifts of the group.} of the
$\infty$-dimensional Lie group ${\rm Ham}(M,\omega)$.

The $L_\infty$ norm on $\mathcal{A}$
$$||F||_\infty={\rm max}\,F-{\rm min}\,F$$
gives rise to the Hofer metric $d$ on ${\rm Ham}(M,\omega)$ in the
following way: we define the Hofer length of a smooth Hamiltonian
path $\alpha : [0,1] \to {\rm Ham}(M,\omega)$ as
$${\rm length}(\alpha):=\int_0^1 ||\dot{\alpha}(t)||_\infty dt=\int_0^1 ||F_t||_\infty dt,$$
where $F_t(x)=F(t,x)$ is the time-dependent Hamiltonian function
generating the path $\alpha$. The Hofer distance $d$ between two
Hamiltonian diffeomorphisms $f$ and $g$ is defined by
$$d(f,g):={\rm inf} \,\{\,{\rm length}(\alpha)\},$$
where the infimum is taken over all Hamiltonian paths $\alpha$
connecting $f$ and $g$.  The Hofer norm $\rho(f)$ is the Hofer
distance between the identity map $id$ and $f$, i.e.
$$\rho(f):=d(id,f).$$

It is easy to check that $d$ is bi-invariant in the sense that
$$d(fh, gh)=d(hf, hg)=d(f, g)$$
for all $f,g,h \in {\rm Ham}(M,\omega)$. The fact that $d$ is
nondegenerate is highly nontrivial. This was proved by Hofer
\cite{Ho} for the case of $\mathbb R^{2n}$, then generalized by
Polterovich \cite{Po} to some larger class of symplectic manifolds,
and finally proved in the full generality by Lalonde and McDuff
\cite{LM} using the following energy-capacity inequality
$$e(S) \geqslant \frac{1}{2} {\rm capacity}(S)$$ for a subset $S$ of
$M$. Here the capacity of $S$ is equal to $\pi r^2$ when $S$ is a
symplectically embedded ball of radius $r$, and is defined in
general as the supremum of the capacities of all symplectically
embedded balls in $S$. The displacement energy $e(S)$ is defined to
be the infimum of the Hofer norms of all $f \in {\rm Ham}(M,\omega)$
such that $f(S) \cap S=\emptyset$.

Note that the energy-capacity inequality provides a lower bound for
the Hofer norm. Namely, we have
$$f(S) \cap S=\emptyset, \,\,{\rm capacity}(S)>c \Longrightarrow \rho(f) >c/2.$$
This fact will be crucial in our proof of Theorem \ref{bicforkt}.

Recall in Definition \ref{bounded} that an element $\phi \in {\rm
Symp}(M,\omega)$ is called unbounded if
$$r(\phi):={\rm sup}\,
\{\,\rho([\phi,f]) \mid f \in {\rm Ham}(M,\omega)\} = \infty.$$ Note
that all Hamiltonian diffeomorphisms are bounded since $r(g)
\leqslant 2 \rho(g)< \infty$ for all $g \in {\rm Ham}(M,\omega)$,
where $\rho(g)$ is the Hofer norm of $g$. According to Proposition
1.2.A in \cite{LP}, $r$ satisfies the triangle inequality $r(\phi
\psi) \leqslant r(\phi)+r(\psi)$. Since ${\rm Ham}(M,\omega)$ is the
kernel of the flux homomorphism, two symplectomorphisms $\phi$ and
$\psi$ have the same flux if and only if they differ by a
Hamiltonian diffeomorphism. Combining these facts, we have the
following

\bigskip

\NI \textbf{Observation A.} \cite{LP} In order to prove
$BI_0(M,\omega)={\rm Ham}(M,\omega)$, it suffices to show that for
each nonzero value $v \in H^1(M, \mathbb R)/\Gamma$, there exists
some unbounded element $\phi \in {\rm Symp}_0(M,\omega)$ with ${\rm
flux}(\phi)=v$.

\bigskip

%%%%%%%%%%%%%%%%%%%%%%%%%%%%%%%%%%%%%%%%%%%%%%%%%%%%%%%%
%%%%%%%%%%%%%%%%%%%%%%%%%%%%%%%%%%%%%%%%%%%%%%%%%%%%%%%
\section{The admissible lift}\label{admissiblelift}

\smallskip

To prove an element $\phi \in {\rm Symp}_0(M,\omega)$ is unbounded,
one has to show that $\rho([\,\phi, f])$ can be arbitrarily large by
choosing different $f \in {\rm Ham}(M,\omega)$. Hence the
energy-capacity inequality will not work directly for closed
manifolds since the capacity of the manifold itself is finite. To go
around this difficulty, we recall the notion of admissible lifts
which was first introduced by Lalonde and Polterovich \cite{LP}. We
shall point out that our definition is slightly different from
theirs, but the two definitions are equivalent.

Let $\pi : (\widetilde{M},\widetilde{\omega}) \to (M,\omega)$ be a
symplectic covering map, i.e. a covering map $\pi$ between two
symplectic manifolds such that $\widetilde{\omega}=\pi^* \omega$.

\begin{definition}
For every $g \in {\rm Ham}(M,\omega)$, assume $g$ is the time-1 map
of the Hamiltonian flow generated by time-dependent Hamiltonian
function $H_t$. An admissible lift $\widetilde{g} \in {\rm
Ham}(\widetilde{M},\widetilde{\omega})$ of $g$ with respect to $\pi$
is defined to be the time-1 map of the Hamiltonian flow generated by
$H_t \circ \pi$.
\end{definition}

\smallskip

\begin{lemma} [existence and uniqueness of admissible lifts]
For all $g \in {\rm Ham}(M,\omega)$, such an admissible lift
$\widetilde{g} \in {\rm Ham}(\widetilde{M},\widetilde{\omega})$
exists and is unique.
\end{lemma}

\begin{proof}
The existence follows from the definition. For the uniqueness, it
suffices to show that the admissible lift $\widetilde{g}$ of $g$ is
independent of the choice of the Hamiltonian function $H_t$.

Note that the choice of $H_t$ is equivalent to the choice of the
Hamiltonian isotopy $g_t$ connecting $id$ to $g$. For every point $p
\in M$, let $$\widetilde{ev}_p : \pi_1({\rm Ham}(M,\omega),id) \to
\pi_1(M,p)$$ be the map induced by the evaluation map $ev_p : {\rm
Ham}(M,\omega) \to M$ which takes $g$ to $g(p)$. It follows from
Floer theory that for all symplectic manifolds $(M,\omega)$, the
induced map $\widetilde{ev}_p$ is trivial, see Chapter 11 \cite{MS1}
for instance. This deep result implies that for any two different
paths $g_t^1$ and $g_t^2$ in ${\rm Ham}(M,\omega)$ connecting $id$
to $g$, $g_t^1(p)$ and $g_t^2(p)$ must be homotopic paths in $M$.
Therefore, for every point $\widetilde{p} \in \widetilde{M}$, the
image $\widetilde{g} (\widetilde{p})$ of $\widetilde{p}$ under
$\widetilde{g}$, being the endpoint of the lift of the path
$g_t(p)$, is independent of the choice of the Hamiltonian isotopy
$g_t$. This proves the uniqueness of admissible lifts.
\end{proof}

\smallskip

For our purposes, we consider the universal cover $\widetilde{M}$ of
$M$. Note that $\widetilde{M}$ is not necessarily compact, and the
admissible lift $\widetilde{g}$ of $g \in {\rm Ham}(M,\omega)$ is
not necessarily compactly supported in $\widetilde{M}$. Instead, it
belongs to ${\rm Ham}_b(\widetilde{M},\widetilde{\omega})$ of
time-$1$ maps of bounded Hamiltonians $\widetilde{M} \times [0, 1]
\to \mathbb R$. The Hofer norm is still well defined and the same
energy-capacity inequality still holds for this setting. This idea
is due to Lalonde and Polterovich \cite{LP}. We shall spell out some
details here for the sake of clarity.

Denote by $(N, \sigma)$ a noncompact symplectic manifold without
boundary. We do not often consider the group ${\rm Ham}(N, \sigma)$
of all Hamiltonian diffeomorphisms with arbitrary support. One
reason in our context is that it would not be possible to define the
Hofer norm on ${\rm Ham}(N, \sigma)$ using the $L_\infty$ norm on
the space $\mathcal A$ of all Hamiltonian functions with arbitrary
support, since not all elements in $\mathcal A$ have finite
$L_\infty$ norms.

One may consider the group ${\rm Ham}_c(N, \sigma)$ of Hamiltonian
diffeomorphisms with compact support. The Hofer norm $\rho$ is well
defined on ${\rm Ham}_c(N, \sigma)$, and the energy-capacity
inequality
$$e_c(S) \geqslant \frac{1}{2} {\rm capacity}(S)$$
is valid as usual, where
$$e_c(S):= {\rm inf}\,\{\, \rho(f)
\mid f \in {\rm Ham}_c(N, \sigma), \, f(S) \cap S=\emptyset\}.$$ As
we have already pointed out, however, this setting is not sufficient
for our purposes since the admissible lift is usually not compactly
supported. Hence we need to consider the larger group ${\rm
Ham}_b(N,\sigma)$ of Hamiltonian diffeomorphisms which are time-$1$
maps of bounded Hamiltonians $H: N \times [0, 1] \to \mathbb R$.
Note that one can not use an arbitrary bounded Hamiltonians $H$,
since the Hamiltonian flow generated by $H$ need not be integrable.
Instead, we only restrict to those bounded Hamiltonians whose flows
are integrable.

The Hofer norm can be defined on ${\rm Ham}_b(N,\sigma)$ exactly the
same way as in Section \ref{hofernorm}. For a subset $S$ of $N$,
define also the bounded displacement energy $e_b(S)$ as
$$e_b(S):= {\rm inf}\,\{\, \rho(f)
\mid f \in {\rm Ham}_b(N, \sigma), \, f(S) \cap S=\emptyset\}.$$
Note that ${\rm Ham}_c(N, \sigma) \subset {\rm Ham}_b(N, \sigma)$
implies $e_b(S) \leqslant e_c(S)$ for any subset $S \subset N$. In
fact, for any compact subset $S$, we have $$e_b(S) = e_c(S).$$ To
prove the other inequality, note that if $f \in {\rm Ham}_b(N,
\sigma)$ displaces a compact subset $S$ from itself, one can easily
construct some cut-off $f_{cut} \in {\rm Ham}_c(N, \sigma)$ of $f$
which still displaces $S$ from itself, and the Hofer norm satisfies
$\rho(f) \geqslant \rho(f_{cut})$. Taking the infimum implies
$e_b(S) \geqslant e_c(S)$.

The above argument implies that the energy-capacity inequality still
holds for the bounded displacement energy. That is
$$e_b(S) \geqslant \frac{1}{2} {\rm capacity}(S).$$
Now back to our discussion about the admissible lift. Note that the
admissible lift $\widetilde{g}$ of $g \in {\rm Ham}(M,\omega)$
belongs to ${\rm Ham}_b(\widetilde{M},\widetilde{\omega})$. And it
follows from the definition of the admissible lift that
$$\rho(g) \geqslant \rho(\widetilde{g})$$
for all $g \in {\rm Ham}(M,\omega)$ and the admissible lift
$\widetilde{g} \in {\rm Ham}_b(\widetilde{M},\widetilde{\omega})$ of
$g$. Here the two $\rho$'s are the Hofer norms on ${\rm
Ham}(M,\omega)$ and ${\rm Ham}_b(\widetilde{M},\widetilde{\omega})$
respectively. Combining the above discussions, we have

\bigskip

\NI \textbf{Observation B.} \cite{LP} To construct $g \in {\rm
Ham}(M,\omega)$ of arbitrarily large Hofer norm, it suffices to make
sure that the unique admissible lift $\widetilde{g} \in {\rm
Ham}_b(\widetilde{M},\widetilde{\omega})$ of $g$ displaces from
itself a symplectic ball in $\widetilde{M}$ of arbitrarily large
capacity.

\bigskip

%%%%%%%%%%%%%%%%%%%%%%%%%%%%%%%%%%%%%%%%%%%%%%%%%%%%%%%%%%%%%%%%%%%
%%%%%%%%%%%%%%%%%%%%%%%%%%%%%%%%%%%%%%%%%%%%%%%%%%%%%%%%%%%%%%%%%%%
\section{Proof of Theorem \ref{bicforkt}}\label{proofbicforkt}

\smallskip

In this section, we prove Theorem \ref{bicforkt}. Recall that $(M,
\omega)$ is the Kodaira-Thurston manifold with the standard
symplectic form $\omega=ds \wedge dt+dx \wedge dy$. Recall also that
$H^1(M, \mathbb R)=\mathbb R\langle ds, dy, dt\rangle $ and the flux
subgroup $\Gamma=\mathbb{Z}\langle ds, dy\rangle $ by Lemma
\ref{cohomology} and Theorem \ref{flux}. In view of Observation A,
to prove $BI_0(M,\omega)={\rm Ham}(M,\omega)$, it suffices to show
that for every nonzero element $v \in H^1(M, \mathbb
R)/\Gamma=\mathbb R/\mathbb{Z} \langle ds, dy\rangle \oplus \mathbb
R\langle dt\rangle $, there exists some unbounded symplectomorphism
with flux equal to $v$. We begin with an explicit construction of
symplectomorphisms with given fluxes.

\begin{lemma}\label{phiabc}
Let $v$ be an element in $H^1(M, \mathbb R)/\Gamma=\mathbb
R/\mathbb{Z} \langle ds, dy\rangle \oplus \mathbb R\langle dt\rangle
$, say $v=\alpha ds+\beta dy+c dt$ where $\alpha, \beta \in \mathbb
R/\mathbb{Z}$ and $c \in \mathbb R$. Then there exists an element
$\phi_{\alpha \beta c} \in {\rm Symp}_0(M,\omega)$ with ${\rm
flux}(\phi_{\alpha \beta c})=v$. Namely,
$$\phi_{\alpha \beta c}(s, t, x, y) = (s+c, t-\alpha, x+\beta, y).$$
\end{lemma}

\begin{proof}
First $\phi_{\alpha \beta c}$ is well-defined. For instance, since
$(s, t, x, y)$ and $(s+1, t, x+y, y)$ represent the same point on
$M$, one has to show that
$$\phi_{\alpha \beta c}(s, t, x, y) \sim  \phi_{\alpha \beta c}(s+1, t, x+y, y).$$
This is true since
$$\phi_{\alpha \beta c} (s, t, x, y)=(s+c, t-\alpha, x+\beta, y),$$
and $$\phi_{\alpha \beta c} (s+1, t, x+y, y)=(s+1+c, t-\alpha,
x+y+\beta, y).$$

\NI It is easy to see that $\phi_{\alpha \beta c}$ preserves
$\omega$, and the obvious isotopy from $id$ to $\phi_{\alpha \beta
c}$ implies that $\phi_{\alpha \beta c} \in {\rm Symp}_0(M,\omega)$.
The calculation for ${\rm flux}(\phi_{\alpha \beta c})=v$ is
straightforward using (\ref{eq:1}) in Section \ref{fluxsubgroup}.
\end{proof}

\smallskip

The following theorem due to Lalonde and Polterovich \cite{LP} is an
important criteria for unbounded symplectomorphisms.

\smallskip

\begin{thm}[Theorem 1.4.A \cite{LP}]\label{displace}
Let $L \subset M$ be a closed Lagrangian submanifold admitting a
Riemannian metric with non-positive sectional curvature, and whose
inclusion in $M$ induces an injection on fundamental groups. Let
$\phi$ be an element in ${\rm Symp}_0(M,\omega)$ such that $\phi(L)
\cap L =\emptyset$. Then $\phi$ is unbounded.
\end{thm}

For the proof, one passes to the universal cover $\widetilde{M}$ of
$M$. The hypothesis implies that the lift of a neighbourhood $U$ of
$L$ has infinite capacity. One then constructs a Hamiltonian isotopy
$f_\tau$ supported in $U$ so that the admissible lift
$\widetilde{[\phi, f_\tau]}$ of the commutator $[\phi, f_\tau]$ will
displace a symplectic ball of arbitrarily large capacity as $\tau$
goes to infinity. This implies $\phi$ is unbounded according to
Observation B. See \cite{LP} for details.

\bigskip

\NI {\textbf{Proof of Theorem \ref{bicforkt}.}} In view of
Observation A, it suffices to show that the symplectomorphisms
$\phi_{\alpha \beta c}$ constructed in Lemma \ref{phiabc} are
unbounded in all cases, as long as the flux $v=\alpha ds+\beta dy+c
dt$ does not vanish. We argue case by case. In the first two cases,
this is a direct consequence of Theorem \ref{displace}.

\NI \textbf{Case 1.} $\alpha \ne 0 \in \mathbb R/ \mathbb{Z}$.

Let $L \subset M$ be the subset of $M$ defined by
$$L := \{(s,t,x,y) \in M \mid t=0,y=0\}.$$
It is easy to check that $L$ is a Lagrangian torus satisfying the
hypothesis of Theorem \ref{displace}, and $\phi_{\alpha \beta c}$
displaces $L$ from itself. Thus $\phi_{\alpha \beta c}$ is
unbounded.

\bigskip

\NI \textbf{Case 2.} $\alpha = 0 \in \mathbb R/ \mathbb{Z}, \, \beta
\ne 0 \in \mathbb R/ \mathbb{Z}$ and $c=0 \in \mathbb R$.

In this case, $\phi_\beta := \phi_{\alpha \beta c}$ maps $(s, t, x,
y)$ to $(s, t, x+\beta, y)$. As in the first case, $\phi_\beta$
displaces from itself a Lagrangian torus $L$ of $M$ defined by
$$L := \{(s, t, x, y) \in M \mid s=0, x=0 \}.$$
We again get $\phi_\beta$ is unbounded in view of Theorem
\ref{displace}.

\bigskip

\NI \textbf{Case 3.} $\alpha = 0 \in \mathbb R/ \mathbb{Z}$ and $c
\ne 0 \in \mathbb R$.

We write $\phi_{\beta c}$ for $\phi_{\alpha \beta c}$ in this case,
$$\phi_{\beta c} := \phi_{\alpha \beta c} : (s, t, x, y) \mapsto (s+c, t, x+\beta, y).$$
Consider two different situations, one of which is simple, while the
other is more complicated.

\bigskip

\NI \textbf{3A.} $\alpha = 0 \in \mathbb R/ \mathbb{Z}$ and $c
\notin \mathbb{Z}$.

As in case 1 and 2, $\phi_{\beta c}$ is unbounded as it displaces
from itself
$$L := \{(s, t, x, y) \in M \mid s=0, x=0 \}.$$

\bigskip

\NI \textbf{3B.} $\alpha = 0 \in \mathbb R/ \mathbb{Z}$ and $c \in
\mathbb{Z} \backslash \{0 \}$.

Note that $(s+c, t, x+\beta, y) \sim (s, t, x +\beta -c y, y)$. So
the map $\phi_{\beta c} : M \to M$ can also be expressed as
$$\phi_{\beta c}(s, t, x, y)=(s, t, x +\beta -c y, y).$$

In contrast to all previous cases where we used the same argument,
here we are facing a difficulty. The trouble is that in this case we
are unable to find a Lagrangian torus of $M$ which is disjoined from
itself by the map $\phi_{\beta c}$. Thus the above argument breaks
down.

To resolve this difficulty, we take $f_\tau$ to be the Hamiltonian
isotopy whose support is in the subset
$$U:=\{(s, t, x, y) \in M \mid |s|<\ep, |x|<\ep \}$$ of $M$.
We require $f_\tau$ to flow only along $y$ and $t$ direction in $U$
and its restriction to
$$V:=\{(s, t, x, y) \in M \mid |s|<\ep/2, |x|<\ep/2 \}$$
is defined by
$$f_\tau(s, t, x, y)=(s, t, x, y-\tau).$$

In the discussion below, $[f, g]:= f g f^{-1} g^{-1}$ stands for the
commutator of $f$ and $g$. Our goal is to show that the unique
admissible lift $\widetilde{[\phi_{\beta c}, f_\tau]}$ of
$[\phi_{\beta c}, f_\tau]$ still displaces from itself a subset of
$\mathbb R^4$ of arbitrarily large capacity when $\tau$ goes to
infinity. For this, we need the following

\begin{lemma}
Let $\phi \in {\rm Symp}_0(M,\omega)$, and $f_\tau$ be a Hamiltonian
isotopy of $M$. Let $\widetilde{\phi} :\widetilde{M} \to
\widetilde{M}$ be any lift of $\phi$, and $\widetilde{[\phi,
f_\tau]}$ and $\widetilde{f}_\tau$ be the unique admissible lift of
$[\phi, f_\tau]$ and $f_\tau$ respectively. Then
$$\widetilde{[\phi, f_\tau]} =
[\widetilde{\phi}, \widetilde{f}_\tau].$$
\end{lemma}

\begin{proof}
Note that $f_\tau$ is Hamiltonian implies $[\phi, f_\tau]$ is
Hamiltonian. So both admissible lifts $\widetilde{[\phi, f_\tau]}$
and $\widetilde{f}_\tau$ make sense. To simplify notation, denote
$$A_\tau:=\widetilde{[\phi, f_\tau]} \,\, \mbox{and} \,\,
B_\tau:=[\widetilde{\phi}, \widetilde{f}_\tau].$$ We want to show
$A_\tau=B_\tau$, which is equivalent to $A_\tau B_\tau^{-1}=id$.
Since $A_\tau$ and $B_\tau$ are both lifts of $[\phi, f_\tau]$,
$A_\tau B_\tau^{-1}$ is the deck transformation of the covering map
$\pi : \widetilde{M} \to M$. Now $A_0 B_0=id$, and $\tau \to A_\tau
B_\tau^{-1}$ is a continuously parametrized path into the discrete
set of all deck transformations. Thus $A_\tau B_\tau^{-1}=id \,$ for
all $\tau$.
\end{proof}

\bigskip

Now back to the proof of Theorem \ref{bicforkt}. To prove
$\phi_{\beta c}$ is unbounded, we need to show that the commutator
$[\phi_{\beta c}, f_\tau]$ has arbitrarily large Hofer norm when
$\tau$ goes to infinity. Let $V_0 \subset \mathbb R^4$ be the subset
of $\mathbb R^4$ defined by
$$V_0:=\{(s,t,x,y) \in \mathbb R^4 \mid |s|< \ep/2, t
\in \mathbb R, |x|<\ep/2, 0<y<\tau/2 \}.$$ Since $V_0$ has
arbitrarily large capacity as $\tau$ goes to infinity, according to
Observation B, it suffices to show that the admissible lift
$\widetilde{[\phi_{\beta c}, f_\tau]}$ of $[\phi_{\beta c}, f_\tau]$
displaces $V_0$ from itself.

For this, denote by $\widetilde{\phi}_{\beta c} : \mathbb R^4 \to
\mathbb R^4$ the preferred lift of the map $\phi_{\beta c}$ such
that
$$\widetilde{\phi}_{\beta c}(s, t, x, y)=(s, t, x +\beta -c y, y).$$
By the above lemma, it suffices to show that
$[\widetilde{\phi}_{\beta c}, \widetilde{f}_\tau] (V_0) \cap V_0 =
\emptyset$, which is equivalent to
$$\widetilde{\phi}_{\beta c}^{-1} \widetilde{f}_\tau^{-1}(V_0) \cap
\widetilde{f}_\tau^{-1} \widetilde{\phi}_{\beta c}^{-1}(V_0)=
\emptyset.$$ Note that the restriction of $\widetilde{f}_\tau$ to
$$\widetilde{V} :=\{(s,t,x,y) \in \mathbb R^4 \mid |s|< \ep/2, t
\in \mathbb R, |x|<\ep/2, y \in \mathbb R \}$$ is defined by
$$\widetilde{f}_\tau(s, t, x, y)=(s, t, x, y-\tau).$$ We have
$$\widetilde{f}_\tau^{-1}(V_0)= \{|s|<
\ep/2, t \in \mathbb R, |x|<\ep/2, \tau<y<3 \tau/2 \}.$$ Hence
$$\widetilde{\phi}_{\beta c}^{-1} \widetilde{f}_\tau^{-1}(V_0)=\{|s|< \ep/2, t \in \mathbb
R, |x +\beta -c y|<\ep/2, \tau<y<3 \tau/2 \}.$$ On the other hand,
$$\widetilde{\phi}_{\beta c}^{-1}(V_0)=\{|s|< \ep/2, t \in
\mathbb R, |x +\beta -c y|<\ep/2, 0<y<\tau/2 \}.$$

\NI Note that in the set $\widetilde{\phi}_{\beta c}^{-1}
\widetilde{f}_\tau^{-1}(V_0)$ we have
$$|x|>|c y| -|\beta| -\ep/2>|c| \tau -|\beta| -\ep/2,$$
and in $\widetilde{\phi}_{\beta c}^{-1} (V_0)$ we have
$$|x|<|c y| +|\beta| +\ep/2<|c| \tau/2 +|\beta| +\ep/2.$$
Thus for sufficiently large $\tau$, these two sets do not share the
same values in $x$ coordinates. Since the flow
$\widetilde{f}_\tau^{-1}$ only changes the $y$ and $t$-coordinates
when restricted to $\widetilde{\phi}_{\beta c}^{-1}(V_0)$, we
conclude
$$\widetilde{\phi}_{\beta c}^{-1} \widetilde{f}_\tau^{-1}(V_0) \cap
\widetilde{f}_\tau^{-1} \widetilde{\phi}_{\beta c}^{-1}(V_0)=
\emptyset.$$ As we have already mentioned above, this implies
$\phi_{\beta c}$ is unbounded in case 3B, which completes the proof
of Theorem \ref{bicforkt}. \qed

\bigskip

\section{Proof of Theorem \ref{bicfortorus}}\label{proofbicfortorus}

\smallskip

We have already mentioned in Section \ref{introduction} that the
bounded isometry conjecture holds for the torus with the standard
symplectic form. In this section we prove Theorem \ref{bicfortorus}
which states that the conjecture holds for the 4-torus with any
linear symplectic form. We begin with a remark on the linear
symplectic form $\omega$ on $\mathbb T^4$.

\begin{rmk}\rm
The 2-form $\omega=\sum_{i < j}\, a_{ij}\, dx_i \wedge dx_j$ on
$\mathbb T^4$ is symplectic, i.e. nondegenerate if and only if
$a_{12}a_{34}-a_{13}a_{24}+a_{14}a_{23} \ne 0$.
\end{rmk}

For each $1 \leqslant i \leqslant 4$, let $\{\phi^i_\theta \} \in
\pi_1({\rm Symp}_0(\mathbb T^4,\omega))$ be the loop of rotations of
$\mathbb T^4$ along $x_i$ direction. Let $\xi_i \in H^1(\mathbb T^4,
\mathbb R)$ be the image of $\{\phi^i_\theta\}$ under the flux
homomorphism. Using (\ref{eq:1}) in Section \ref{fluxsubgroup}, one
easily gets
$$\xi_i := {\rm flux}(\{\phi^i_\theta\})=\displaystyle \sum_{j=1}^4
\,a_{ij} \, dx_j.$$ Here we take the convention that
$a_{ij}=-a_{ji}$. In particular, $a_{ii}=0$.

\begin{lemma}
For the $4$-torus with the linear symplectic form $\omega:=\sum_{i <
j}\, a_{ij} dx_i \wedge dx_j$, the flux subgroup $\Gamma \subset
H^1(\mathbb T^4, \mathbb R)$ is generated by the above $\xi_i's$
over $\mathbb Z$. That is, $\Gamma = \mathbb Z \langle \xi_1, \xi_2,
\xi_3, \xi_4 \rangle$.
\end{lemma}

\begin{proof}
According to Lemma \ref{commute}, we have the following commutative
diagram for the manifold $(\mathbb T^4, \omega)$.

$$\begin{CD}
  \pi_1({\rm Symp}_0(\mathbb T^4,\omega)) @>\widetilde{ev}_s>>H_1(\mathbb T^4, \mathbb Z)    @>{\rm PD}>>H^3(\mathbb T^4, \mathbb Z) \\
  @VidVV             &&                      @VV\centerdot {\rm vol}(\mathbb T^4)V\\
  \pi_1({\rm Symp}_0(\mathbb T^4,\omega)) @>{\rm flux}>>H^1(\mathbb T^4, \mathbb R) @>\wedge[\omega]>>H^3(\mathbb T^4, \mathbb R).
\end{CD}$$

\bigskip

\NI Note that $\widetilde{ev}_s$ is surjective, and $\wedge[\omega]
: H^1(\mathbb T^4, \mathbb R) \to H^3(\mathbb T^4, \mathbb R)$ is an
isomorphism. Note also that ${\rm vol}(\mathbb
T^4)=a_{12}a_{34}-a_{13}a_{24}+a_{14}a_{23}$. It follows from a
similar argument as in the proof of Theorem \ref{flux} that $\xi_i
(1 \leqslant i \leqslant 4)$ span the flux subgroup $\Gamma$ over
$\mathbb Z$.
\end{proof}

Now let $\phi \in Symp_0(\mathbb T^4,\omega)$ such that
$$\phi(x_1, x_2, x_3, x_4)=(x_1+\alpha_1, x_2+\alpha_2, x_3+\alpha_3,
x_4+\alpha_4)$$ where $\alpha_i \in \mathbb R/\mathbb Z$ for $1
\leqslant i \leqslant 4$. Then
$${\rm flux}(\phi)= \displaystyle \sum_{i=1}^4 \, \alpha_i \,
\xi_i.$$

Recall that in view of Observation A in Section \ref{hofernorm}, to
prove Theorem \ref{bicfortorus}, it suffices to show $\phi$ is
unbounded as long as at least one $\alpha_i \in \mathbb R/\mathbb Z$
is nonzero. One may attempt to apply Theorem \ref{displace} by
showing $\phi$ disjoins some Lagrangian torus $L \subset \mathbb
T^4$ from itself. For a general symplectic form $\omega$, however,
there may not exist any such Lagrangian torus in $\mathbb T^4$.
Nevertheless, we can still prove $\phi$ is unbounded using the
following

\begin{lemma}\label{nocontractible}
Let $(M,\omega)$ be an aspherical symplectic manifold. Let $f_\tau
\in {\rm Ham}(M, \omega)$ be the flow generated by an autonomous
Hamiltonian which has no nonconstant contractible orbits. Then the
Hofer norm $\rho(f_\tau)$ goes to infinity as $\tau$ goes to
infinity.
\end{lemma}

This result can be found in Oh \cite{Oh}, Schwarz \cite{Sch} and
Kerman-Lalonde \cite{KL}. The main idea of the argument is that the
Hofer norm is bounded from below by the spectral norm, while the
spectral norm of such $f_\tau$ grows linearly with respect to
$\tau$.

\bigskip

\NI {\bf Proof of Theorem \ref{bicfortorus}.} Let $\phi \in
Symp_0(\mathbb T^4,\omega)$ such that
$$\phi(x_1, x_2, x_3, x_4)=(x_1+\alpha_1, x_2+\alpha_2, x_3+\alpha_3,
x_4+\alpha_4).$$ As discussed above, it suffices to show $\phi$ is
unbounded when at least one $\alpha_i \in \mathbb R/\mathbb Z$ is
nonzero. Assume $\alpha_1 \ne 0$ without loss of generality. Thus
$\phi (U) \cap U = \emptyset$ where $U \subset \mathbb T^4$ is
defined by
$$U := \{(x_1, x_2, x_3, x_4) \in \mathbb T^4 \mid |x_1|<
\epsilon \}.$$ for sufficiently small $\ep$.

Let $H$ be a time-independent Hamiltonian function of $\mathbb T^4$
supported in $U$. Denote by $f_\tau$ the (autonomous) Hamiltonian
flow generated by $H$. Since $\phi (U) \cap U = \emptyset$, we know
that $[\phi, f_\tau]:=\phi f_\tau \phi^{-1} f_\tau^{-1}$ is also an
autonomous Hamiltonian flow supported in the union of two disjoint
sets $U \cup \phi (U)$. If we further require that $H$ depend only
on the first coordinate $x_1$, using the fact that $\omega$ is a
linear symplectic form, we conclude that $[\phi, f_\tau]$ has no
nonconstant contractible orbits. Thus it follows from Lemma
\ref{nocontractible} that the Hofer norm $\rho([\phi, f_\tau])$ goes
to infinity as $\tau$ goes to infinity. Hence $\phi$ is unbounded in
the sense of Definition \ref{bounded}. \qed

\bigskip

%%%%%%%%%%%%%%%%%%%%%%%%%%%%%%%%%%%%%%%%%%%%%%%%%%%%%%%%%%%%%%%%%%%%%%%
%%%%%%%%%%%%%%%%%%%%%%%%%%%%%%%%%%%%%%%%%%%%%%%%%%%%%%%%%%%%%%%%%%%%%%%
\section{The Kodaira-Thurston manifold with linear symplectic
forms}\label{sec8}

\smallskip

So far we have studied bounded isometries for the Kodaira-Thurston
manifold with the standard symplectic form and for the 4-torus with
all linear symplectic forms. In particular, we have shown that the
bounded isometry conjecture holds in both cases. In this section we
will study the same question for the Kodaira-Thurston manifold with
all linear symplectic forms.

\begin{question}\label{bicforktlinear}
Does the bounded isometry conjecture hold for the Kodaira-Thurston
manifold with all linear symplectic forms?
\end{question}

We expect the answer to be positive. Although we are not able to
give a complete proof yet at this time, we shall provide some
partial results below. We begin by describing the linear symplectic
forms on the Kodaira-Thurston manifold $M$. Recall that it follows
from Lemma \ref{cohomology} that $H^2(M, \mathbb R)$ is of rank $4$,
generated by $\gamma \wedge ds, \, \gamma \wedge dy, \, ds \wedge
dt, \, \mbox{and} \,dy \wedge dt$ where $\gamma=dx-s dy$. We
consider linear 2-forms
$$\omega_{a b e f}:=a \gamma \wedge ds + b\gamma \wedge d y + e ds \wedge dt + f dy \wedge dt.$$
Note that $\omega_{a b e f}$ is a symplectic form if and only if $b
e - a f \ne 0$. In particular, the standard symplectic form
corresponds to $b = e =1$ and $ a = f =0$. The following lemma on
the flux subgroup generalizes Theorem \ref{flux}.

\begin{lemma}\label{fluxlinear}
The flux subgroup $\Gamma \subset H^1(M,\mathbb R)$ of the
Kodaira-Thurston manifold with the linear symplectic form $\omega_{a
b e f}$ has rank $2$ over $\mathbb{Z}$. More precisely, we have
$\Gamma=\mathbb{Z}\langle e ds + f dy, \,a ds + b dy \rangle $.
\end{lemma}

\begin{proof}
The proof follows the same lines as that of Theorem \ref{flux}.
According to Lemma \ref{commute}, we have the following commutative
diagram.

$$\begin{CD}
  \pi_1({\rm Symp}_0(M,\omega_{abef})) @>\widetilde{ev}_s>>H_1(M, \mathbb Z)    @>{\rm PD}>>H^3(M, \mathbb Z) \\
  @VidVV             &&                      @VV\centerdot {\rm vol}(M)V\\
  \pi_1({\rm Symp}_0(M,\omega_{abef})) @>{\rm flux}>>H^1(M, \mathbb R) @>\wedge[\omega_{abef}]>> H^3(M, \mathbb R).
\end{CD}$$

\bigskip

As in the proof of Theorem \ref{flux}, the image of
$\widetilde{ev}_s$ in $H_1(M, \mathbb Z)$ is contained in
$\mathbb{Z}\langle\frac{\partial}{\partial t}\rangle$. Note that
$PD(\frac{\partial}{\partial t})=-\gamma \wedge dy \wedge ds$, where
$\gamma=dx-sdy$. Now look at the map $\wedge\omega_{abef} : H^1(M,
\mathbb R) \to H^3(M, \mathbb R)$,
\begin{equation}\nonumber
 \begin{aligned}
   ds \mapsto ds \wedge \omega_{abef} &=b \gamma \wedge dy \wedge ds - f dy \wedge ds \wedge dt =b \gamma \wedge dy \wedge ds,\\
     dt \mapsto dt \wedge \omega_{abef} &=a \gamma \wedge ds \wedge dt + b \gamma \wedge dy \wedge dt,\\
       dy \mapsto dy \wedge \omega_{abef} &= -a \gamma \wedge dy \wedge ds + e dy \wedge ds \wedge dt=-a \gamma \wedge dy \wedge ds.
        \end{aligned}.
          \end{equation}

\smallskip

\NI Here we have used the fact that the 3-form $dy \wedge ds \wedge
dt = d (\gamma \wedge dt)$ is exact, so it vanishes on the
cohomology level. Since ${\rm vol}(M)=be-af \ne 0$, we conclude by
tracing the diagram that the flux subgroup $\Gamma \subset H^1(M,
\mathbb R)$ is contained in $\mathbb{Z}\langle e ds + f dy, a ds+ b
dy\rangle $. Note that the fact $be-af \ne 0$ implies that $e ds + f
dy$ and $a ds+ b dy$ are linearly independent. An explicit
construction shows that $\Gamma$ is actually equal to
$\mathbb{Z}\langle e ds + f dy, a ds + b dy \rangle$. Namely, we
take two elements $\{\phi_\theta\}$ and $\{\psi_\theta\}$ in
$\pi_1({\rm Symp}_0(M,\omega_{abef}))$ such that
$$\phi_\theta(s, t, x, y)=(s, t-\theta, x, y), 0 \leqslant\theta \leqslant1, $$
$$\psi_\theta(s, t, x, y)=(s, t, x+\theta, y), 0 \leqslant\theta \leqslant1.$$
A straightforward calculation using (\ref{eq:1}) in Section
\ref{fluxsubgroup} shows that ${\rm flux}(\{\phi_\theta \}) =e ds +
f dy$ and ${\rm flux}(\{\psi_\theta \})=a ds + b dy$.
\end{proof}

As in Lemma \ref{phiabc}, we explicitly construct below
symplectomorphisms with given fluxes.

\begin{lemma}\label{phiabclinear}
Let $v$ be an element in $$H^1(M, \mathbb R)/\Gamma=\mathbb
R/\mathbb{Z} \langle e ds + f dy, \, a ds + b dy \rangle \oplus
\mathbb R\langle dt\rangle,$$ say $$v=\alpha (e ds + f dy) +\beta (a
ds + b dy)+c (be - af) dt$$ where $\alpha, \beta \in \mathbb
R/\mathbb{Z}$ and $c \in \mathbb R$. Then there exists $\phi_{\alpha
\beta c} \in {\rm Symp}_0(M,\omega_{abef})$ with ${\rm
flux}(\phi_{\alpha \beta c})=v$. Namely,
$$\phi_{\alpha \beta c}(s, t, x, y) = (s+ b c, t-\alpha, x+\beta - a c s, y - a c).$$
\end{lemma}

\begin{proof}
First $\phi_{\alpha \beta c}$ is well-defined. For instance, since
$(s, t, x, y)$ and $(s+1, t, x+y, y)$ represent the same point in
$M$, one has to show that
$$\phi_{\alpha \beta c}(s, t, x, y) \sim  \phi_{\alpha \beta c}(s+1, t, x+y, y).$$
This is true since
$$\phi_{\alpha \beta c}(s, t, x, y) = (s+ b c, t-\alpha, x+\beta - a c s, y - a c)$$
and $$\phi_{\alpha \beta c} (s+1, t, x+y, y)=(s+ 1 + b c, t-\alpha,
x+ y + \beta - a c (s+1), y - a c)$$ also represent the same point.
One can check that $\phi_{\alpha \beta c}^* \omega_{abef} =
\omega_{abef}$, and the obvious isotopy from $id$ to $\phi_{\alpha
\beta c}$ implies that $\phi_{\alpha \beta c} \in {\rm
Symp}_0(M,\omega_{abef})$.

It remains to show that ${\rm flux}(\phi_{\alpha \beta c})=v$. Note
that $\phi_{\alpha \beta c}$ is the time-1 map of the flow generated
by the time-independent symplectic vector field
$$X:= bc \frac{\partial}{\partial s}- \alpha
\frac{\partial}{\partial t} + (\beta - acs) \frac{\partial}{\partial
x} - ac \frac{\partial}{\partial y}.$$ Using (\ref{eq:1}) in Section
\ref{fluxsubgroup}, we have
\begin{equation}\nonumber
\begin{split}
{\rm flux} (\phi_{\alpha \beta c})&=\iota(X)\,\omega_{a b e f}\\
 &= \iota (bc
\frac{\partial}{\partial s}- \alpha \frac{\partial}{\partial t}
+ (\beta - acs) \frac{\partial}{\partial x} - ac \frac{\partial}{\partial y})\, \omega_{a b e f}\\
 &= -abc (dx - s dy) + bce dt + \alpha e ds + \alpha f dy \\
 &+ a (\beta - acs) ds + b (\beta - acs) dy + a^2 c sds + abc dx - acf dt\\
 &= \alpha (e ds + f dy) + \beta (a ds + b dy) + c(be - af) dt\\
 &= v.
\end{split}
 \end{equation}
\end{proof}

To answer Question \ref{bicforktlinear}, one has to check whether
$\phi_{\alpha \beta c}$ constructed in Lemma \ref{phiabclinear} is
always unbounded whenever its flux $v$ is nonzero in $H^1(M, \mathbb
R)/\Gamma$. This is in general a very hard question. In the
remaining of this section, we will give a proof for some known
cases. For the unknown cases, we will try to point out what
difficulty is involved.

\bigskip

\NI {\bf Case 1}: $\alpha \ne 0 \in \mathbb R/\mathbb{Z}$. In this
case we will prove $\phi_{\alpha \beta c}$ is always unbounded. Note
that $\phi_{\alpha \beta c}(U) \cap U = \emptyset$ where $U \subset
M$ is defined by
$$U := \{(s, t, x, y) \in M \mid |t|<
\epsilon \}$$ for sufficiently small $\ep$. We will apply Lemma
\ref{nocontractible} as in the proof of Theorem \ref{bicfortorus}.
Recall that the only thing we need to do is to construct
time-independent Hamiltonian $H$ supported in $U$ whose flow has no
nonconstant contractible orbits. This follows from a tedious but
straightforward calculation which asserts that
$$\iota\,(X)\,\omega_{abef}=dt$$ where
$$X:=\frac{1}{be-af}(-as\frac{\partial}{\partial x}-a
\frac{\partial}{\partial y}+b \frac{\partial}{\partial s}).$$

\NI Note that this is actually a special case of the construction in
Lemma \ref{phiabclinear}. And the fact that $X$ is a well defined
vector field on $M$ follows from the equivalence relation $(s,t,x,y)
\sim (s+1,t,x+y,y)$. Since $a$ and $b$ can not be both zero, if we
further require $H$ to depend only on the $t$-coordinates, we know
that the Hamiltonian flow generated by $H$ will have no nonconstant
contractible orbits. Therefore $\phi_{\alpha \beta c}$ is always
unbounded in this case.

\bigskip

\NI {\bf Case 2}: $\alpha = 0 \in \mathbb R/\mathbb{Z}$ and $c \ne 0
\in \mathbb R$. First we assume $ac$ and $bc$ are not both integers.
Note that this is always the case when the ratio $a : b$ is
irrational. Under this assumption, $\phi_{\beta c}:=\phi_{\alpha
\beta c}$ is unbounded in view of Theorem \ref{displace} as it
disjoins a Lagrangian torus
$$L:=\{(s,t,x, y) \in M \mid s =0, \, y=0 \}.$$

If the ratio $a : b$ is rational, then there exists $c \ne 0$ such
that both $ac$ and $bc$ are integers. In this case, using the
equivalence relation $(s,t,x,y) \sim (s+1,t,x+y,y)$, we can write
the map
$$\phi_{\beta c}:(s, t, x, y) \mapsto (s+ b c, t, x+\beta - a c s, y - a c)$$
as
$$\phi_{\beta c}:(s, t, x, y) \mapsto (s, t, x+\beta - a c s-b c y, y).$$
It is natural to attempt the admissible lift argument as in Case 3B
of Theorem \ref{bicforkt} for the standard Kodaira-Thurston
manifold. One would try to construct a Hamiltonian isotopy
$\widetilde{f}_\tau$ on $\mathbb R^4$ supported in
$$\widetilde{U}:=\{(s, t, x, y) \in \mathbb R^4 \mid |es+fy|<\ep, |x|<\ep \}$$
which flows only along $s$ and $y$ directions, and whose restriction
to
$$\widetilde{V}:=\{(s, t, x, y) \in \mathbb R^4 \mid |es+fy|<\ep/2, |x|<\ep/2 \}$$
is defined by
$$\widetilde{f}_\tau(s, t, x, y)=(s+f \tau, t, x, y-e \tau).$$
Note that the above construction allows us to show that the lift
$$\widetilde{\phi}_{\beta c} :(s,
t, x, y) \mapsto (s, t, x+\beta - a c s-b c y, y)$$ of $\phi_{\beta
c}$ is unbounded on the universal cover level. For this, one would
argue as in Case 3B of Theorem \ref{bicforkt}, that the commutator
$[\widetilde{\phi}_{\beta c}, \widetilde{f}_\tau]$ displaces some
subset $V_0 \subset \mathbb R^4$ of arbitrarily large capacity with
respect to the symplectic form
$\widetilde{\omega}_{abef}:=\pi^*\omega_{abef}$. Namely,
$$V_0:=\{|es+fy|<\ep/2, \, t
\in \mathbb R,\, |x|<\ep/2, 0<as+by<|be-af|\tau/2 \}.$$ The problem
here is that $\widetilde{f}_\tau$ does not descend to a Hamiltonian
isotopy on $M$. Note that in proving $\phi_{\beta c}$ itself is
unbounded, it is crucial to have such a Hamiltonian isotopy on $M$,
not just on the universal cover $\mathbb R^4$. Hence this case is
still unsolved.

\bigskip

\NI {\bf Case 3}: $\alpha = 0 \in \mathbb R/\mathbb{Z}$, $c = 0 \in
\mathbb R$ and $\beta \ne 0 \in \mathbb R/\mathbb{Z}$. In this case,
the map $\phi_\beta:= \phi_{\alpha \beta c}$ has the simple form
$$\phi_\beta :(s, t, x, y) \mapsto (s, t, x+\beta, y).$$

\NI We do not know in general how to prove $\phi_\beta$ is unbounded
for this seemingly easy case. The difficulty in applying Theorem
\ref{displace} is that the obvious torus
$$L := \{(s, t, x, y) \in M \mid s=0, x=0 \}$$ displaced by
$\phi_\beta$ is not necessarily Lagrangian with respect to all
symplectic forms $\omega_{abef}$. If we assume $f=0$, then $L$ is
actually a Lagrangian torus, and $\phi_\beta$ will be unbounded in
view of Theorem \ref{displace}.

Note also that Lemma \ref{nocontractible} does not work here either
since our situation here is different from Case 1 above. The main
reason is that
$$U := \{(s, t, x, y) \in M \mid |x|< \epsilon \}$$
is not a well defined set in $M$. Thus one can no longer apply Lemma
\ref{nocontractible} by constructing a time-independent Hamiltonian
$H$ supported in $U$ whose flow has no nonconstant contractible
orbits.

\end{document}